\documentclass[11pt]{article}
\usepackage{verbatim}
\usepackage{amsfonts,amsmath,amsthm}
\def\LaTeX{L\kern -.36em\raise .3ex\hbox{\sc a}\kern -.15em T\kern -.1667em%
\lower .7ex\hbox{E}\kern -.125em X}

\usepackage{amssymb}
\usepackage{float}
\usepackage{epsfig}
\usepackage{amsmath}
\usepackage[english]{babel}
\usepackage{a4}
\usepackage{amstext}
\usepackage{mathrsfs}
\usepackage{amsthm}
\usepackage{fancyhdr}
\usepackage{latexsym}
\usepackage{amsmath}
\usepackage{amsfonts}
\usepackage{amssymb}
\usepackage{fancyhdr}

\usepackage[usenames]{color}
\definecolor{red}{rgb}{1.0,0.0,0.0}

\definecolor{blu}{rgb}{0.0,0.0,1.0}

\def\to{\longrightarrow}

\def\norm{{\| \kern -.05em | }}

\newtheorem{theorem}{Theorem}[section]
\newtheorem{proposition}[theorem]{Proposition}
\newtheorem{lemma}[theorem]{Lemma}

\theoremstyle{definition}
\newtheorem{definition}[theorem]{Definition}
\newtheorem{remark}[theorem]{Remark}

\newtheorem{hypothesis}[theorem]{Hypothesis}

\author{Salvatore Federico \and Ben Goldys \and Fausto Gozzi}

\begin{document}

\title{
HJB Equations for the Optimal Control of Differential Equations with
Delays and State Constraints, I: Regularity and
Applications. \footnote{This work was partially supported by an
Australian Research Council Discovery Project}}
\author{
Salvatore Federico\footnote{\thinspace Salvatore Federico,
Dipartimento di Scienze Economiche ed Aziendali, Facolt\`{a} di
Economia, Libera Universit\`{a} internazionale degli studi sociali
``Guido Carli'', viale Romania 32, 00197 Roma, Italy. Email:
\texttt{sfederico at luiss.it}.} \and Ben Goldys\footnote{\thinspace
Ben Goldys, Schoolo of Mathematics and Statistics, University of New
South Wales, Sydney, Australia. Email: \texttt{B.Goldys at
unsw.edu.au}.} \and Fausto Gozzi\footnote{\thinspace Fausto Gozzi
(corresponding author), Dipartimento di Scienze Economiche ed
Aziendali, Facolt\`{a} di Economia, Libera Universit\`{a}
internazionale degli studi sociali ``Guido Carli'', viale Romania
32, 00197 Roma, Italy. Email: \texttt{fgozzi at luiss.it}.} }

\maketitle
\begin{abstract}
We study a class of optimal control problems with state constraints
where the state equation is a differential equation with delays.
This class includes some problems arising in economics, in
particular the so-called models with time to build, see
\cite{AseaZak,BambiJEDC,KydlandPrescott}. We embed the problem in a
suitable Hilbert space $H$ and consider the associated
Hamilton-Jacobi-Bellman (HJB) equation. This kind of
infinite-dimensional HJB equation has not been previously studied
and is difficult due to the presence of state constraints and the
lack of smoothing properties of the state equation. Our main result
on the regularity of solutions to such a HJB equation seems to be
completely new. More precisely we prove that the value function is
continuous in a sufficiently big open set of $H$, that it solves in
the viscosity sense the associated HJB equation and it has
continuous classical derivative in the direction of the ``present''.
This regularity result is the starting point to define a feedback
map in classical sense, which gives rise to a candidate optimal
feedback strategy for the problem. The study of verification
theorems and of the closed loop equation will be the subject of the
forthcoming paper \cite{FGG2}.
\end{abstract}


\bigskip

\noindent \textbf{A.M.S. Subject Classification}: 34K35,49L25,
49K25.

\newpage

\tableofcontents

\section{Introduction}

The main purpose of this paper is to prove a $C^1$ regularity result
for a class of first order infinite dimensional HJB equations
associated to the optimal control of deterministic delay equations
arising in economic models.

The $C^1$ regularity of solutions of the HJB equations arising
in deterministic optimal control theory is a crucial issue to solve
in a satisfactory way the control problems. Indeed, even in finite
dimension, in order to obtain the optimal strategies in feedback form one usually needs
the existence of an appropriately defined gradient of the solution. It is possible to prove
verification theorems and representation of optimal feedbacks in the
framework of viscosity solutions, even if the gradient is not
defined in classical sense (see e.g. \cite{BCD,YongZhou}), but
this is usually not satisfactory in applied problems since the closed loop
equation becomes very hard to treat in such cases.

The need of $C^1$ regularity results for HJB equations is
particularly important in infinite dimension since in this case
verification theorems in the framework of viscosity solutions are
rather weak and in any case not applicable to problems with state
constraints (see e.g \cite{FGSJCA, LiYong}). To the best of our
knowledge $C^1$ regularity results for first order HJB equation have
been proved by method of convex regularization introduced by Barbu
and Da Prato \cite{BDPbook} and then developed by various authors
(see e.g. \cite{BDP1,BDP2,BDPP,BP,DB1,DB2,FagAMO,GSICON,GVorau}). All
these results do not hold in the case of state constraints and, even
without state constraints, do not cover problems where the state
equation is a nonlinear differential equation with delays.  In the papers
\cite{CDB1,CDB2,FagSICON} a class of state constraints problems is
treated using the method of convex regularization but the $C^1$ type
regularity is not proved.

In this paper we deal with a class of optimal control problems
where, given a control $c\ge 0$ the  state $x$ satisfies the
following delay equation
$$
\begin{cases}
x'(t)=rx(t)+f_0\left(x(t),\int_{-T}^0a(\xi)x(t+\xi)d\xi\right)-c(t),\\
x(0)=\eta_0, \ x(s)=\eta_1(s), \ s\in[-T,0),
\end{cases}
$$
with state constraint $x(\cdot)>0$.
The objective is to maximize the functional
$$J(\eta; c(\cdot)):=\int_0^{+\infty} e^{-\rho t}
\left[U_1(c(t))+ U_2(x(t))\right]\,dt, \ \ \ \rho> 0,$$ over the set
of the admissible controls $c$. We will give a more precise
formulation of the problem and of the hypotheses on the data in the
next section. For the moment we observe that this kind of problems
arises in various economic models (see e.g.
\cite{AseaZak,BambiJEDC,KydlandPrescott}), where the authors study
optimal growth in presence of time-to-build (i.e. delay in the
production due to the need of time to build new products) and cannot
be treated with the existing theory except for very special cases
(see the three papers just quoted).

Using a standard approach (see e.g. \cite{BenDap}) we embed the
problem in an infinite dimensional control problem in the Hilbert
space $H=\mathbb{R}\times L^2([-T,0];\mathbb{R})$, where intuitively
speaking $\mathbb R$ describes the ``present'' and
$L^2([-T,0];\mathbb{R})$ describes the ``past'' of the system. The
associated Hamilton-Jacobi-Bellman (HJB) equation in $H$ has not
been previously studied and is difficult due to the presence of
state constraints and the lack of smoothing properties of the state
equation.

We prove that the value function is continuous in a sufficiently big
open set of $H$ (Proposition \ref{A^-1}), that it solves in the
viscosity sense the associated HJB equation (Theorem \ref{TH:visc})
and it has continuous classical derivative in the direction of the
``present'' (Theorem \ref{TH:reg}). This regularity result is enough
to define the formal optimal feedback strategy in classical sense,
since the objective functional only depends on the ``present''.
The method we use to prove regularity is completely different from
the one of convex regularization mentioned above. Indeed, it is based
on a finite dimensional result of Cannarsa and Soner
\cite{CannarsaSoner} (see also \cite{BCD}, pag. 80) that exploits
the concavity of the data and the strict convexity of the
Hamiltonian to prove the continuous differentiability of the
viscosity solution of the HJB equation. Generalizing to the infinite
dimensional case such result is not trivial as the definition of
viscosity solution in this case strongly depends (via  the unbounded
differential operator $A$ contained in the state equation) on the
structure of the problem. In particular we need to establish
specific properties of superdifferential that are given in
Subsection \ref{sec:super}.

We believe that such a method could be also used to analyze other
problems featuring concavity of the data and strict convexity of the
Hamiltonian.

We finally observe that, even with our regularity result at hand the
Verification Theorems and the study of the Closed Loop Equation
associated to the problem are not trivial at all and will be the
subject of the forthcoming paper \cite{FGG2}.

The plan of the paper is as follows. Section \ref{Sec:pr} is devoted
to set up the problem in DDE form giving main assumptions and some
preliminary results (in Subsection \ref{sec:value}) that are proved
directly without using the infinite dimensional setting. In Section
3 we rewrite the problem in the infinite dimensional setting and
prove existence and uniqueness of solutions of the state equation
(Subsection \ref{subs:mild}), continuity of the value function
(Subsection \ref{sec:continuity}) and useful properties of
superdifferentials (Subsection \ref{sec:super}). In Section 4 we
apply the dynamic programming in the infinite dimensional context to
get our main results: we prove that the value function is a
viscosity solution of the HJB equation (Subsection \ref{sec:visc})
and then we prove a regularity result for viscosity solutions of HJB
(Subsection \ref{sec:reg}).

\section{Setup of the control problem and preliminary results}\label{Sec:pr}
In this section we will formally define the control delay problem
giving some possible applications of it. We will use the notations
\[L^2_{-T}:=L^2([-T,0];\mathbb{R}),\quad\mathrm{and}\quad
W^{1,2}_{-T}:=W^{1,2}([-T,0];\mathbb{R}).\]
We will denote by $H$ the Hilbert space
$$H:=\mathbb{R}\times L^2_{-T},$$ endowed with the inner product
$$\langle\cdot,\cdot\rangle= \langle\cdot,\cdot\rangle_{\mathbb{R}}+
\langle\cdot,\cdot\rangle_{L^2_{-T}},$$
and the norm
\[\|\cdot\|^2=|\cdot|^2_{\mathbb R}+\|\cdot\|^2_{L^2_{-T}}.\]
We will denote by
$\eta:=(\eta_0,\eta_1(\cdot))$ the generic element of this space.
For convenience we set also$$H_+:=(0,+\infty)\times L^2_{-T}, \ \ \ \ \
H_{++}:=(0,+\infty)\times\{\eta_1(\cdot)\in L^2_{-T} \ | \
\eta_1(\cdot)\geq 0 \ a.e.\}.$$
\begin{remark}\label{rm:meglioH+}
Economic motivations we are mainly interested in (see
\cite{AseaZak,BambiJEDC,KydlandPrescott} and Remark \ref{fin} above)
require to study the optimal control problem with the initial
condition in $H_{++}$. However, the set $H_{++}$ is not convenient
to work with, since its interior with respect to the
$\|\cdot\|$-norm is empty. That is why we enlarge the problem and
allow the initial state belonging to the class $H_+$.
\hfill$\blacksquare$
\end{remark}

For $\eta\in H_+$, we consider an optimal control of the following
differential delay equation:
\begin{equation}\label{eqstate}
\begin{cases}
x'(t)=rx(t)+f_0\left(x(t),\int_{-T}^0a(\xi)x(t+\xi)d\xi\right)-c(t),\\
x(0)=\eta_0, \ x(s)=\eta_1(s), \ s\in[-T,0),
\end{cases}
\end{equation}
with state constraint $x(\cdot)>0$ and control constraint $c(\cdot)\geq 0$.
We set up the following assumptions on the functions $a,f_0$.
\begin{hypothesis}\label{f_0,a}
\begin{itemize}
\item[]
\item $a(\cdot)\in W^{1,2}_{-T}$ is such that $a(\cdot)\geq 0$ and $a(-T)=0$;
\item $f_0:[0,\infty)\times \mathbb{R}\rightarrow \mathbb{R}$
 is jointly concave, {nondecreasing} with respect to the second variable, Lipschitz continuous with Lipschitz constant $C_{f_0}$, and
 \begin{equation}\label{f_0strict}
 f_0(0,y)>0, \ \ \ \ \forall y>0.
 \end{equation}
\end{itemize}
\hfill $\blacksquare$
\end{hypothesis}

\begin{remark}\label{rm:Bambipossibile}
In the papers \cite{AseaZak,BambiJEDC,KydlandPrescott} the pointwise
delay is used. We cannot treat exactly this case for technical
reasons that are explained in Remark \ref{rm:noBambi} below. However
we have the freedom of choosing the function $a$ in a wide class and
this allows to take account of various economic phenomena. Moreover
we can approximate the pointwise delay with suitable sequence of
functions $\{a_n\}$ getting convergence of the value function and
constructing $\varepsilon$-optimal strategies (this approximation
procedure is the object of the forthcoming paper \cite{FGG2}).
\hfill$\blacksquare$
\end{remark}

From now on we will assume that $f_0$ is extended to a Lipschitz
continuous map on $\mathbb{R}^2$ setting
$$f_0(x,y):=f_0(0,y), \ \ \ \mbox{for} \ x<0.$$
For technical reasons, which will be clear in Subsection
\ref{sec:continuity}, we work with the case $r>0$, noting that
nevertheless the case $r\leq 0$ can be treated by shifting the
linear part of the state equation. Indeed in this case we can
rewrite the state equation taking for example as new coefficient for
the linear part $\tilde{r}=1$ and shifting the nonlinear term
defining $\tilde{f}_0(x,y)=f_0(x,y)-(1-r) x$.
\par\medskip\noindent
We say that a {function} $x:[-T,\infty)\to\mathbb R^+$ is a solution to equation \eqref{eqstate} if $x(t)=\eta_1(t)$ for $t\in[-T,0)$ and
\[x(t)=\eta_0+\int_0^trx(s)ds+\int_0^tf_0\left(x(s),\int_{-T}^0a(\xi)x(s+\xi)d\xi\right)ds-\int_0^tc(s)ds,\quad t\ge 0.\]
\begin{theorem}\label{daprato}
 For any given $\eta \in H_+$, $c(\cdot)\in L^1_{loc}([0,+\infty);\mathbb{R^+})$,  equation \eqref{eqstate} admits a unique solution that is absolutely continuous {on $[0,+\infty)$}.
\end{theorem}
\textbf{Proof.} Let $K=\sup_{\xi\in[-T,0]}a(\xi)$. For any $t\geq 0$, $z^1, z^2\in C([-T,t];\mathbb{R})$, we have

\begin{multline*}
\int_0^t\left[r|z_1(s)-z_2(s)|+\left|f_0\left(z_1(s),\int_{-T}^0a(\xi)z_1(s+\xi)\right)-f_0\left(z_2(s),\int_{-T}^0a(\xi)z_2(s+\xi)\right)\right|\right]ds\\
\leq \int_0^t \left[r|z_1(s)-z_2(s)|+ C_{f_0}\left[ |z_1(s)-z_2(s)|+K\int_{-T}^0|z_1(s+\xi)-z_2(s+\xi)|d\xi\right]\right]ds\\
\leq \int_0^t \left[(r+C_{f_0})|z_1(s)-z_2(s)|+ C_{f_0}K\int_{-T}^t|z_1(\xi)-z_2(\xi)|d\xi\right]ds\\
\leq  (r+C_{f_0})\int_0^t|z_1(s)-z_2(s)|ds+ tC_{f_0}K\int_{-T}^t|z_1(\xi)-z_2(\xi)|d\xi\\
\leq [(r+C_{f_0})+tC_{f_0}K]\int_{-T}^t|z_1(\xi)-z_2(\xi)|d\xi\\
\leq [(r+C_{f_0})+tC_{f_0}K](t+T)^{1/2}\left(\int_{-T}^t|z_1(\xi)-z_2(\xi)|^2d\xi\right)^{1/2}.
\end{multline*}
Thus the claim follows by  Theorem 3.2, pag. 246,  of \cite{BenDap}. \hfill$\square$\\

We denote by $x(\cdot; \eta, c(\cdot))$ the unique solution of
(\ref{eqstate}) with initial point $\eta\in H_+$ and under the
control $c(\cdot)$. We emphasize that this solution actually
satisfies pointwise only the integral equation associated with
\eqref{eqstate}; it satisfies \eqref{eqstate} in differential form
only for almost every $t\in[0,+\infty)$. \vspace{.3cm}

For $\eta\in H_+$ we define the class of the admissible controls
starting from $\eta$ as $$\mathcal{C}(\eta):= \{ c(\cdot)\in
L^1_{loc}([0,+\infty);\mathbb{R^+}) \ | \ x(\cdot;  \eta, c(\cdot))
> 0\}.$$ Setting $x(\cdot):=x(\cdot,\,; \eta,c(\cdot))$, the problem
consists in maximizing the functional
$$J(\eta; c(\cdot)):=\int_0^{+\infty} e^{-\rho t} \left[U_1(c(t))+ U_2(x(t))\right]\,dt, \ \ \ \rho> 0,$$
over the set of the admissible strategies. \vspace{.3cm}

The following will be standing assumptions on the utility functions
$U_1$, $U_2$, holding throughout the whole paper.
\begin{hypothesis}\label{hyp:utility}
\begin{itemize}
\item[]
\item[(i)] $U_1\in C([0,+\infty);\mathbb{R})\cap C^2 ((0,+\infty);\mathbb{R})$,  $U_1'>0$, $U_1'(0^+)=+\infty$, $U_1''<0$ and $U_1$ is bounded.
\item[(ii)] $U_2\in C((0,+\infty);\mathbb{R})$ is increasing, concave, bounded from above. Moreover
\begin{equation}\label{ipotesiU_2}
\int_0^{+\infty} e^{-\rho t} U_2\left(e^{-C_{f_0} t}\right)dt>-\infty.
\end{equation}
\end{itemize}
\vspace{-0,50truecm}\hfill $\blacksquare$
\end{hypothesis}
Since $U_1$, $U_2$ are bounded from above, the previous functional is well-defined for any $\eta\in H_+$ and $c(\cdot)\in\mathcal{C}(\eta)$.
We set
$$\bar{U}_1:=\lim_{s\rightarrow+\infty} U_1(s), \ \ \ \bar{U}_2:=\lim_{s\rightarrow+\infty} U_2(s).$$

\begin{remark}\label{utility}
We give some comments on Hypothesis \ref{hyp:utility} and on the
structure of the utility in the objective functional.
\begin{enumerate}
\item {Through the whole paper the case $U_2\equiv 0$ is allowed.  Therefore, the case of an objective  functional
depending only on consumption (as in \cite{AseaZak,BambiJEDC,KydlandPrescott}) is allowed. }
\item
The assumption that $U_1$, $U_2$ are bounded from above is done for
simplicity to avoid too many technicalities. It guarantees that the
value function is bounded from above and this fact simplifies
arguments in the following parts of this paper. Similarly the
assumption that $U_1$ is bounded from below guarantees that the
value function is bounded from below. We think that it is possible
to replace such assumptions with more general conditions relating
the growth of $U_1$, $U_2$, the value of $\rho$ and the parameters
of the state equation. Typically, such a condition requires that
$\rho$ is sufficiently large.
\item
All utility functions bounded from below satisfy \eqref{ipotesiU_2}.
Also $U_2(x)=\log(x), \ \ \ U_2(x)=x^{\gamma}, \
\gamma>-\frac{\rho}{C_{f_0}},$ satisfy \eqref{ipotesiU_2}. Note also
that \eqref{ipotesiU_2}  is equivalent to
\begin{equation*}
\int_0^{+\infty} e^{-\rho t} U_2\left(\xi e^{-C_{f_0}
t}\right)dt>-\infty, \ \ \ \forall \xi>0.
\end{equation*}
\item When $r<0$, then in \eqref{ipotesiU_2} we have to replace $C_{f_0}$ with $|r|+C_{f_0}$.
\item If we assume that
\begin{equation}\label{conditionaltra}
\exists \delta>0 \ \mbox{such that} \ rx+f_0(x,0)\geq 0, \forall x\in (0,\delta],
\end{equation}
then the assumption \eqref{ipotesiU_2} can be suppressed.
Since  \eqref{f_0strict} implies $f_0(0,0)\geq 0$, \eqref{conditionaltra} occurs for example if $x\mapsto rx+f_0(x,0)$ is nondecreasing, therefore in particular if $r\geq 0$ and $f_0$ depends only on the second variable.
\item We think that it should be possible to treat also the case of a utility function $U$ depending on the couple $(c,x)$; it should be enough to replace Hyphotesis \ref{hyp:utility} with the following assumptions:
\begin{itemize}
\item $U:[0,+\infty)\times(0,+\infty)\rightarrow \mathbb{R}$ is concave and increasing with respect to both the variables and bounded from above;
\item for any $x>0$ the function $U(\cdot,x)$ belongs to the class $C([0,+\infty);\mathbb{R})\cap C^2 ((0,+\infty);\mathbb{R})$, $U'(\cdot, x)>0$, $U'(0^+,x)=+\infty$, $U''(\cdot,x)<0.$
\end{itemize}
\end{enumerate}
\end{remark}

\begin{remark}\label{fin}
The control problem described above covers also the following optimal consumption problem. We may think of the dynamics defined in (\ref{eqstate}) as the dynamics of the bank account driven by a contract which takes into account the past history of the accumulation of wealth. Such a situation arises when the bank offers to the customer an interest rate $r$ smaller than the market spot rate $r^M$ and as a compensation, it provides a premium on the past of the wealth (this may happen e.g. for pension funds). Then the following equation is a possible simple model of the evolution of the bank account under such a contract:
\begin{equation*}
\begin{cases}
x'(t)=rx(t)+g_0\left(\int_{-T}^0 a(\xi)x(t+\xi)d\xi\right)-c(t),\\
x(0)=\eta_0, \ x(s)=\eta_1(s), \ s\in[-T,0),
\end{cases}
\end{equation*}
where $g_0:\mathbb{R}\rightarrow \mathbb{R}$ is a concave, Lipschitz continuous and strictly increasing function such that $g_0(0)\geq 0$. Dependence on the past is an incentive for the customer to keep his investments with the bank for a longer period of time in order to receive gains produced by the term $g_0$. Here we  assume the point of view of the customer and we think it is interesting to study the behaviour of the optimal consumption in this setting, comparing it with the one coming from the classical case, which corresponds to set $r=r^M$, $g\equiv 0$.

We think also that our technique could be adapted to cover optimal advertising  model with nonlinear memory effects (see e.g. \cite{GoMaSa} on this subject in a stochastic framework).
\hfill$\blacksquare$
\end{remark}

 For $\eta\in H_+$ the value function of our problem is defined by
\begin{equation}\label{value}
V(\eta):=\sup_{c(\cdot)\in\mathcal{C}(\eta)} J(\eta, c(\cdot)),
\end{equation}
with the convention $\sup \emptyset =-\infty.$
 The domain of the value function is the set
$$\mathcal{D}(V):=\{\eta\in H_+ \ | \ V(\eta)>-\infty\}.$$
Due to the assumptions on $U_1$, $U_2$ we directly get that $V\leq \frac{1}{\rho}(\bar{U}_1+\bar{U}_2)$.

\subsection{Preliminary results}\label{sec:value}
In this subsection we investigate some first qualitative properties
of the delay state equation and of the value function.
\begin{lemma}[Comparison]\label{comparison}
 Let $\eta\in H_+$ and let $c(\cdot)\in L^1_{loc}([0,+\infty);\mathbb{R^+})$. Let $x(t)$, $t\geq 0$, be an absolutely continuous function satisfying almost everywhere the differential inequality
\begin{equation*}
\begin{cases}
x'(t)\leq rx(t)+f_0\left(x(t),\int_{-T}^0a(\xi)x(t+\xi)d\xi\right)-c(t),\\
x(0)\leq\eta_0, \ x(s)\leq\eta_1(s), \ \mbox{for a.e.} \ s\in[-T,0).
\end{cases}
\end{equation*}
Then
$x(\cdot)\leq x(\cdot;{\eta},{c}(\cdot)).$\\

\end{lemma}
\textbf{Proof.}
Set $\bar{a}:=\sup_{\xi\in[-T,0]}|a(\xi)|$, $y(\cdot):=x(\cdot;\eta,c(\cdot))$ and $h(\cdot):=[x(\cdot)-y(\cdot)]^+$. We must show that $h(\cdot)=0$. Let  $\varepsilon>0$ be fixed and such that $\varepsilon C_{f_0}\bar{a} T  e^{\varepsilon (r+C_{f_0})}\leq 1/2$ and let $M:=\max_{t\in[0,\varepsilon]} h(t)$.
By monotonicity of $f_0$
we get
\begin{equation}\label{casc}
f_0\left(x(t),\int_{-T}^0a(\xi)x(t+\xi)d\xi\right)\leq f_0\left(x(t),\int_{-T}^0a(\xi)y(t+\xi)d\xi+ \bar{a}TM\right), \ \ \ \ \mbox{for} \ t\in[0,\varepsilon].
\end{equation}
Define, for $n\in\mathbb{N}$,
$$\varphi_n(x):=\begin{cases}
0, \ \ \ \ \ \ \ \ \ \ \ \ \ \mbox{for}  \ x\leq 0,\\
nx^2, \  \ \ \ \ \ \ \ \ \ \mbox{for} \ x\in(0,1/2n],\\
x-1/4n, \ \ \ \mbox{for} \ x>1/2n.
\end{cases}$$
The sequence $(\varphi_n)_{n\in\mathbb{N}}\subset C^1(\mathbb{R};\mathbb{R})$ is such that
$$
\begin{cases}
\varphi_n(x)=\varphi_n'(x)= 0,  \ \ \mbox{for every} \ x\in (-\infty,0], \ n\in\mathbb{N},\\
0\leq \varphi'_n(x)\leq 1, \ \ \mbox{for every} \ x\in\mathbb{R}, \ n\in\mathbb{N},\\
\varphi_n(x)\rightarrow x^+, \ \ \mbox{ uniformly on} \ x\in\mathbb{R},\\
\varphi'_n(x)\rightarrow 1, \ \mbox{for} \ x\in(0,+\infty).
\end{cases}
$$
We can write for $t\in[0,\varepsilon]$, taking into account \eqref{casc},
\begin{eqnarray*}
\varphi_n(x(t)-y(t))&=&\varphi_n (x(0)-\eta_0)+\int_{0}^t \varphi'_n(x(s)-y(s))[x'(s)-y'(s)]ds\\
&\leq &  \int_{0}^t \varphi'_n(x(s)-y(s))\Bigg[r(x(s)-y(s))\\
&& \ \ \ +f_0\left(x(s),\int_{-T}^0a(\xi)x(s+\xi)d\xi\right)-f_0\left(y(s),\int_{-T}^0a(\xi)y(s+\xi)d\xi\right)\Bigg]ds\\
&\leq &  \int_{0}^t \varphi'_n(x(s)-y(s))\Bigg[r(x(s)-y(s))\\
&& \ \ \ +f_0\left(x(s),\int_{-T}^0a(\xi)y(s+\xi)d\xi+\bar{a} TM\right)-f_0\left(y(s),\int_{-T}^0a(\xi)y(s+\xi)d\xi \right)\Bigg]ds\\
&\leq &  \int_{0}^t \varphi'_n(x(s)-y(s))\Big[(r+C_{f_0})|x(s)-y(s)|+C\bar{a} TM\Big]ds.
\end{eqnarray*}
Letting $n\rightarrow \infty$ we get
$$
h(t)\leq \int_{0}^t (r+C_{f_0})\,h(s)ds + C_{f_0}\bar{a} TM t\leq \int_{0}^t (r+C_{f_0})\,h(s)ds + C_{f_0}\bar{a} TM \varepsilon.
$$
Therefore by  Gronwall's Lemma we get
$$h(t)\leq \varepsilon C_{f_0}\bar{a}T M e^{\varepsilon (r+C_{f_0})}, \ \ \ \mbox{for} \ t\in[0,\varepsilon],$$
and by definition of $\varepsilon$
$$h(t)\leq \frac{M}{2}, \ \ \ \mbox{for} \ t\in[0,\varepsilon].$$
This shows that $M=0$, i.e. that $h=0$ on  $[0,\varepsilon]$.
Iterating the argument, since $\varepsilon$ is fixed, we get $h\equiv0$ on $[0,+\infty)$, i.e. the claim.
\hfill$\square$



\begin{proposition}\label{piu}
We have
$$H_{++}\subset\mathcal{D}(V), \ \ \ \ \ \mathcal{D}(V)=\{\eta\in H_+ \ | \ 0\in \mathcal{C}(\eta)\}.$$
\end{proposition}
\textbf{Proof.}
Let $\eta\in H_{++}$ and set $x(\cdot):=x(\cdot;\eta,0)$. By assumption, $x(0)=\eta_0>0$ and until $x(t)>0$ we have
$$x'(t)= rx(t)+f_0\left(x(t),\int_{-T}^0a(\xi)x(t+\xi)d\xi\right)\geq rx(t)+f_0(x(t),0).$$
Since $f_0(0,0)\geq 0$ and $f_0(\cdot,0)$ is Lipschitz continuous (with Lipschitz constant $C_{f_0}$),  we get
$$x'(t)\geq -C_{f_0} \,x(t), \ \ \ \mbox{ until} \ x(t)>0.$$
This fact forces to be  $$\inf\{t\geq 0 \ | \ x(t)=0\}=+\infty,$$
and $x(t)\geq \eta_0e^{-C_{f_0} t}$ for any $t\geq 0$, which proves the inclusion $H_{++}\subset\mathcal{D}(V)$ thanks to \eqref{ipotesiU_2}.

Now let $\eta\in \mathcal{D}(V)$; then, by definition of $\mathcal{D}(V)$, there exists $c(\cdot)\in\mathcal{C}(\eta)$. By Lemma \ref{comparison}  $0\in\mathcal{C}(\eta)$, so that we have the inclusion $\mathcal{D}(V)\subset \{\eta\in H_+ \ | \ 0\in\mathcal{C}(\eta)\}$. Conversely let $\eta\in H_+$ be such that $0\in\mathcal{C}(\eta)$. Then, by definition of $\mathcal{C}(\eta)$, we have $\inf_{t\in[0,T]} x(t;\eta,0)\geq \xi>0$. Repeating the argument used above, we get $x(t;\eta,0)\geq \xi e^{-C_{f_0}(t-T)}$ for $t\geq T$, so that $\eta\in\mathcal{D}(V)$ and the proof is  complete. \hfill$\square$

\begin{remark}
It is straightforward to see that the proof of Proposition \ref{piu} above works if we replace the assumption \eqref{ipotesiU_2} with the assumption \eqref{conditionaltra}.
\hfill $\blacksquare$
\end{remark}
\begin{definition}
(i)
Let $\eta\in \mathcal{D}(V)$. An admissible control $c^*(\cdot)\in\mathcal{C}(\eta)$ is said to be optimal for the initial state $\eta$ if $J(\eta; c^*(\cdot))=V(\eta)$. In this case the corresponding state trajectory $x^*(\cdot):=x(\cdot; \eta, c^*(\cdot))$ is said to be an optimal trajectory and the couple $(x^*(\cdot),c^*(\cdot))$ is said an optimal couple.

(ii)  Let $\eta\in \mathcal{D}(V), \ \varepsilon >0$; an admissible control $c^\varepsilon(\cdot)\in\mathcal{C}(\eta)$ is said $\varepsilon$-optimal for the initial state $\eta$ if $J(\eta; c^\varepsilon(\cdot))>V(\eta)-\varepsilon$.  In this case the corresponding state trajectory $x^\varepsilon(\cdot):=x(\cdot; \eta, c^\varepsilon(\cdot))$ is said an $\varepsilon$-optimal trajectory and the couple $(x^\varepsilon(\cdot),c^\varepsilon(\cdot))$ is said an $\varepsilon$-optimal couple.   \hfill$\square$
\end{definition}

\begin{proposition}  The set $\mathcal{D}(V)$ is convex and
the value function $V$  is  concave on $\mathcal{D}(V)$.
\end{proposition}
\textbf{Proof.} Let $\eta,\bar{\eta}\in \mathcal{D}(V)$ and set, for $\lambda\in[0,1]$,  $\eta_\lambda= \lambda \eta+(1-\lambda)\bar{\eta}$. For $\varepsilon>0$, let $c^\varepsilon(\cdot)\in\mathcal{C}(\eta)$, $\bar{c}^\varepsilon(\cdot)\in\mathcal{C}(\bar{\eta})$ be two controls $\varepsilon$-optimal for the initial states $\eta,\bar{\eta}$ respectively. Set $x(\cdot)=x(\cdot, \eta, c^\varepsilon(\cdot))$, $\bar{x}(\cdot)=x(\cdot; \eta, \bar{c}^\varepsilon(\cdot))$,  $c^\lambda(\cdot)=\lambda c^\varepsilon(\cdot)+(1-\lambda)\bar{c}^\varepsilon(\cdot)$. Finally set $x_\lambda (\cdot)=\lambda x(\cdot)+(1-\lambda)\bar{x}(\cdot)$. Let us write the dynamics for $x_\lambda(\cdot)$:
\begin{eqnarray*}
x'_\lambda(t)&=&\lambda x'(t)+(1-\lambda)\bar{x}'(t)\\
&=& \lambda \left[rx(t)+f_0\left(x(t), \int_{-T}^0a(\xi)x(t+\xi)d\xi\right)-c^\varepsilon (t)\right]\\
&&+ (1-\lambda) \left[r\bar{x}(t)+f_0\left(\bar{x}(t), \int_{-T}^0a(\xi)\bar{x}(t+\xi)d\xi\right)-\bar{c}^\varepsilon (t)\right]\\
&\leq &rx_\lambda(t)+f_0\left(x_\lambda(t), \int_{-T}^0a(\xi)x_\lambda(t+\xi)d\xi\right)-c^\lambda (t),
\end{eqnarray*}
where the  inequality follows from the concavity of $f_0$
with initial condition $\eta_\lambda$. Let $x(\cdot;\eta_\lambda, c^\lambda(\cdot))$ be a solution of the equation
\[x^\prime(t)=rx(t)+f_0\left(x(t),\int_{-T}^0a(\xi)x(t+\xi)d\xi\right)-c^\lambda(t).\]
Since $x_\lambda(\cdot)> 0$ by construction, by Lemma \ref{comparison} we have  $x(\cdot;\eta_\lambda, c^\lambda(\cdot))\geq x_\lambda(\cdot)> 0$. This shows that $c^\lambda(\cdot)\in \mathcal{C}(\eta_\lambda)$. By  concavity of $U_1$, $U_2$ and by monotonicity of $U_2$ we get

$$V(\eta_\lambda)\geq J(\eta_\lambda; c^\lambda(\cdot))\geq \lambda  J(\eta; c^\varepsilon(\cdot)) + (1-\lambda)J(\eta; \bar{c}^\varepsilon(\cdot))> \lambda V(\eta)+(1-\lambda) V(\bar{\eta})-\varepsilon.$$
Since $\varepsilon$ is arbitrary, we get the claim.\hfill$\square$ \\

By assumptions of monotonicity of the utility functions and by Lemma \ref{comparison} we obtain the following result.
\begin{proposition}\label{cresc}
The function $\eta\mapsto V(\eta)$ is nondecreasing in the sense that
$$\eta_0\geq \bar{\eta}_0, \ \eta_1(\cdot)\geq \bar{\eta}_1(\cdot)\Longrightarrow V(\eta_0,\eta_1(\cdot))\geq V(\bar{\eta}_0,\bar{\eta}_1(\cdot)).$$
\hfill$\square$
\end{proposition}


Indeed, the value function is strictly increasing in the first variable:
\begin{proposition}\label{strict}
We have the following statements:
\begin{enumerate}
\item $V(\eta)< \frac{1}{\rho}(\bar{U}_1+\bar{U}_2)$ for any $\eta\in H_+$.
\item  $\lim_{\eta_0\rightarrow +\infty} V(\eta_0,\eta_1(\cdot))=\frac{1}{\rho}(\bar{U}_1+\bar{U}_2)$, for all $\eta_1(\cdot)\in L^2_{-T}$.
\item $V$ is strictly increasing with respect to the first variable.
\end{enumerate}
\end{proposition}
\textbf{Proof. 1.}
Let $\eta\in \mathcal{D}(V)$ and set
$$\bar{a}:=\sup_{\xi\in[-T,0]} a(\xi), \ \ \ \ p:=\sup_{\xi\in [0,T]} x(\xi; \eta,0), \ \ \ \ q:=\int_{-T}^0\eta_1^+(\xi)d\xi.$$
Let $c(\cdot)\in\mathcal{C}(\eta)$ and set $x(\cdot):=x(\cdot; \eta, c(\cdot))$; by comparison criterion we have $x(t)\leq p$ in $[0,T]$.

{Since $f_0$ is Lipschitz continuous, there exists $C>0$ such that $f_0(x,y)\leq C(1+|x|+|y|)$ for all $x,y\in \mathbb{R}$. Therefore,
for $t\in[0,T]$, we can write, considering the state equation in integral form,
$$x(t)\leq \eta_0+ rTp+ T C (1+|p|+|\bar{a}(Tp+q)|)-\int_0^T c(\tau)d\tau.$$
Call $K:=\eta_0+ rTp+ T C(1+|p|+|\bar{a}(Tp+q)|)$;}
since $c(\cdot)\in \mathcal{C}(\eta)$, we have $x(t)> 0$ in $[0,T]$, so that
$$\int_0^T c(\tau)d\tau \leq K.$$
Denoting by $m$ the Lebesgue measure, this means that
$$m\{\tau\in[0,T] \ | \ c(\tau)\leq 2K/T\}\geq T/2.$$
Therefore (in the next inequality, since $e^{-\rho t}$ is decreasing, we suppose without loss of generality that $c(\cdot)\leq 2K/T$ on $\left[\frac{T}{2},T\right]$)
\begin{eqnarray*}
\int_0^{+\infty} e^{-\rho t} U_1(c(t))dt& \leq & \int_0^{T/2} e^{-\rho t} U_1(c(t))dt+\int_{T/2}^{T} e^{-\rho t} U_1(2K/T) dt+\int_{T}^{+\infty} e^{-\rho t} U_1(c(t))dt\\
&\leq &\frac{\bar{U}_1}{\rho}-\int_{T/2}^T e^{-\rho t} \Big(\bar{U}_1- U_1(2K/T)\Big)dt.
\end{eqnarray*}
Since the quantity $\bar{U}_1-  U_1(2K/T)$ is strictly positive and does not depend on $c(\cdot)$, the claim is proved.

\vspace{.1cm}
\textbf{2.} {For given  ${\eta}_1(\cdot)\in L^2_{-T}$, let $K>0$, $M>0$ and let us define the control
$$
c(t):=\begin{cases}
M, \ \ \ \mbox{if} \ t\in [0,K],\\
0, \  \ \ \ \ \mbox{if} \ t>K.
\end{cases}
$$
Take $\eta_0>0$.
Since $f_0$ is Lipschitz continuous and nondecreasing with respect to the second variable, we can see that, until it is positive, $x(t; (\eta_0,\eta_1(\cdot)), c(\cdot))$ satisfies the differential inequality
$$
\begin{cases}
x'(t)\geq  -C(1+x(t) +q)- M,\\
x(0)=\eta_0,
\end{cases}
$$
for some $C>0$,
where $$q:=\left(\sup_{\xi\in [-T,0]}a(\xi)\right) \left(\int_{-T}^0\eta_1^-(\xi)d\xi\right).$$}
This actually shows that, for any $M>0$, $K>0$, $R>0$, we can find $\eta_0$ such that $c(\cdot)\in\mathcal{C}(\eta_0,\eta_1(\cdot))$  and $x(\cdot; (\eta_0,\eta_1(\cdot)),c(\cdot))\geq R$ on $[0,K]$. By the arbitrariness of $M,K,R$ the claim is proved.

\vspace{.1cm}
\textbf{3.} Fix $\eta_1(\cdot)$; we know that $\eta_0\mapsto V(\eta_0,\eta_1(\cdot))$ is concave and increasing. If it is not strictly increasing, then it has to be constant on an half line $[k, +\infty)$, but this contradicts the first two claims. \hfill$\square$\\

\section{The delay problem rephrased in  infinite dimension}\label{sec:infinite}
Our aim is to apply the dynamic programming technique in order to solve the control problem described in the previous section. However, this approach requires a markovian setting. That is why we will reformulate the problem as an infinite-dimensional control problem. Let $\hat{n}=(1,0)\in H_+$ and
let us consider, for $\eta\in H$ and $c(\cdot)\in L^1([0,+\infty);\mathbb{R}^+)$, the following evolution equation in the space $H$:
\begin{equation}\label{infinitestate}
\begin{cases}
X'(t)=AX(t)+ F(X(t))-c(t)\hat{n},\\
X(0)=\eta\in H_+.
\end{cases}
\end{equation}
In the equation above:
\begin{itemize}
\item $A:\mathcal{D}(A)\subset H\longrightarrow H$ is an unbounded operator defined by $A(\eta_0, \eta_1(\cdot)):= (r\eta_0, \eta_1'(\cdot))$ on
$$\mathcal{D}(A):=\{ \eta\in H \ | \ \eta_1(\cdot)\in W^{1,2}_{-T}, \ \eta_1(0)=\eta_0\};$$
\item $F:H\longrightarrow H$ is a Lipschitz continuous map defined by
\begin{equation*}
F(\eta_0,\eta_1(\cdot)):=\left(f\left(\eta_0,\eta_1(\cdot)\right),0\right),
\end{equation*}
where $f(\eta_0,\eta_1(\cdot)):=f_0\left(\eta_0,\int_{-T}^0a(\xi)\eta_1(\xi)d\xi\right)$.
\end{itemize}
It is well known that $A$ is the infinitesimal generator of a
strongly continuous semigroup $(S(t))_{t\geq 0}$ on $H$; its explicit expression is given by
$$
S(t)(\eta_0,\eta_1(\cdot))=\left(\eta_0e^{rt},I_{[-T,0]}(t+\cdot) \
\eta_1(t+\cdot)+ I_{[0,+\infty)}(t+\cdot) \ \eta_0e^{r(t+\cdot)}\right);
$$
About the estimate on the norm of the semigroup, we have
\begin{eqnarray*}
\|S(t)\eta\|^2&\leq&\left|\eta_0e^{rt}\right|^2+2\int_{-T}^0\left|
I_{[-T,0]}(t+\zeta) \
\eta_1(t+\zeta)\right|^2d\zeta\\&&+2\int_{-T}^0\left|
I_{[0,+\infty)}(t+\zeta) \ \eta_0e^{r(t+\zeta)}\right|^2d\zeta\\
&\leq& ((3+2T))e^{2rt}\|\eta\|^2,
\end{eqnarray*}
i.e.
\begin{equation}\label{semigroup}
\|S(t)\|_{\mathcal{L}(H)}\leq Me^{\omega t},
\end{equation}
where $M=(3+2T)$, $\omega=2r$.
\subsection{Mild solutions of the state equation}\label{subs:mild}
In this subsection we give a definition of the \emph{mild solution}
to  (\ref{infinitestate}),  prove the existence and uniqueness of
such a solution and the equivalence between the one dimensional
delay problem and the infinite dimensional one.
\begin{definition}
A mild solution of
(\ref{infinitestate})
is a function $X\in C([0,+\infty);H)$ which satisfies the integral equation
\begin{equation}\label{mildcap3}
X(t)=S(t)\eta+\int_0^t S(t-\tau)F(X(\tau))d\tau
+\int_0^t c(\tau)S(t-\tau)\hat{n}\,d\tau.
\end{equation}
\end{definition}
\begin{theorem}\label{existmild}
For any $\eta\in H$, there exists a unique mild solution of (\ref{infinitestate}).
\end{theorem}
\textbf{Proof.}
Due to the Lipschitz continuity of $F$ and to \eqref{semigroup},
the proof is the usual standard application of the fixed point theorem (see e.g. \cite{BenDap}).\hfill$\square$\\

We denote by $X(\cdot;\eta, c(\cdot))=(X_0(\cdot;\eta,c(\cdot)), X_1(\cdot;\eta,c(\cdot)))$ the unique solution to (\ref{infinitestate}) for the initial state $\eta\in H$ and under the control $c(\cdot)\in L^1([0,+\infty);\mathbb{R}^+)$.
The following equivalence result justifies our approach.
\begin{proposition} \label{link}
Let $\eta\in H_+$, $c(\cdot)\in\mathcal{C}(\eta)$ and let $x(\cdot)$, $X(\cdot)$ be respectively the unique solution to (\ref{eqstate}) and the unique mild solution to (\ref{infinitestate}) starting from $\eta$ and under the control $c(\cdot)$. Then, for any $t\geq 0$, we have the equality in $H$  $$X(t)=\left(x(t),x(t+\xi)_{\xi\in [-T,0]}\right).$$
\end{proposition}
\textbf{Proof.} Let $x(\cdot)$ be a solution of \eqref{eqstate} and let $Z(\cdot):=(x(\cdot),x(\cdot+\zeta)|_{\zeta\in[-T,0]})$. Then $Z(\cdot)$ belongs to the space $C([0,+\infty);H)$ because the function $[0,+\infty)\ni t\mapsto x(t)\in\mathbb{R}$ is (absolutely) continuous. Therefore,
it remains to prove that
$Z(t)=(Z_0(t),Z_1(t))$ satisfies
(\ref{infinitestate}) and then the claim will follow by uniqueness. For the first component we have to
verify that, for any $t\geq 0$,
$$
Z_0(t)=e^{rt}\eta_0+\int_0^te^{r(t-\tau)}f_0\left(Z_0(\tau), \int_{-T}^0 a(\xi)Z_1(\tau)(\xi)d\xi\right)d\tau-\int_0^te^{r(t-\tau)}c(\tau)d\tau,
$$
i.e. that
$$
x(t)=e^{rt}\eta_0+\int_0^te^{r(t-\tau)}f_0\left(x(\tau), \int_{-T}^0 a(\xi)x(\tau+\xi)d\xi\right)d\tau-\int_0^te^{r(t-\tau)}c(\tau)d\tau,
$$
but this follow from the assumption that $x(\cdot)$
is a solution to (\ref{eqstate}).

\noindent
For the second component, taking into account that
$I_{[0,+\infty)}(t+\cdot-\tau)=I_{[\tau,+\infty)}(t+~\cdot)$, we
have to verify, for any $t\geq 0$, for a.e. $\zeta\in[-T,0]$,
\begin{eqnarray*}
Z_1(t)(\zeta)&=&I_{[-T,0]}(t+\zeta)\eta_1(t+\zeta)+I_{[0,+\infty)}(t+\zeta)\eta_0e^{r(t+\zeta)}\\
&&+\int_0^tI_{[\tau,+\infty)}(t+\zeta)\,\,e^{r(t+\zeta-\tau)}f_0\left(Z_0(\tau), \int_{-T}^0 a(\xi)Z_1(\tau)(\xi)d\xi\right)d\tau\\
&&-\int_0^t I_{[\tau,+\infty)}(t+\zeta)\,\, e^{r(t+\zeta-\tau)}c(\tau)d\tau,\\
\end{eqnarray*}
i.e., for any $t\geq 0$, for a.e. $\zeta\in[-T,0]$,
\begin{eqnarray}\label{opop2}
x(t+\zeta)&=&I_{[-T,0]}(t+\zeta)\eta_1(t+\zeta)+I_{[0,+\infty)}(t+\zeta)\,\,\eta_0e^{r(t+\zeta)}\nonumber\\
&&+\int_0^tI_{[\tau,+\infty)}(t+\zeta)\,\,e^{r(t+\zeta-\tau)}f_0\left(x(\tau), \int_{-T}^0 a(\xi)x(\tau+\xi)d\xi\right)d\tau\\
&&-\int_0^t I_{[\tau,+\infty)}(t+\zeta)\,\, e^{r(t+\zeta-\tau)}c(\tau)d\tau.\nonumber
\end{eqnarray}
For
$\zeta\in[-T,0]$ such that $t+\zeta\in[-T,0]$, (\ref{opop2}) reduces to
$$x(t+\zeta)=\eta_1(t+\zeta)$$
and this is true since $\eta_1$ is the initial condition of (\ref{eqstate}). If
$\zeta\in[-T,0]$ is such that $t+\zeta\geq 0$, then (\ref{opop2}) reduces to
$$
x(t+\zeta)=\eta_0e^{r(t+\zeta)}
+\int_0^{t+\zeta} e^{r(t+\zeta-\tau)}f_0\left(x(\tau), \int_{-T}^0 a(\xi)x(\tau+\xi)d\xi\right)d\tau
-\int_0^{t+\zeta} e^{r(t+\zeta-\tau)}c(\tau)d\tau.
$$
Setting $u:=t+\zeta$ this equality becomes, for $u\geq 0$,
$$
x(u)=x_0e^{ru}+\int_0^t e^{r(u-\tau)}f_0\left(x(\tau), \int_{-T}^0 a(\xi)x(\tau+\xi)d\xi\right)d\tau
-\int_0^t e^{r(u-\tau)}c(\tau)d\tau.
$$
Again this is true because $x(\cdot)$ solves (\ref{eqstate}).
\hfill$\square$
\\

\subsection{Continuity of the value function}\label{sec:continuity}
In this subsection we prove a continuity property of the value function that will be useful to investigate the geometry of its superdifferential in the next subsection. To this end we recall that the generator $A$ of the semigroup $(S(t))_{t\geq 0}$ has bounded inverse in $H$ given by
\[A^{-1}\left(\eta_0,\eta_1\right)(s)=\left(\frac{\eta_0}{r},\frac{\eta_0}{r}-\int_s^0\eta_1(\xi)d\xi\right),\quad s\in [-T,0].\]
It is well known that $A^{-1}$ is compact in $H$. It is also clear that $A^{-1}$ is an isomorphism of $H$ onto $\mathcal{D}(A)$ endowed with the graph norm.
\par\noindent
We define the $\|\cdot\|_{-1}$-norm on $H$ by $$\|\eta\|_{-1}:=\|A^{-1}\eta\|.$$
In the next proposition we characterize the adjoint operator $A^*$ and its domain $\mathcal{D}(A^*)$.
\begin{proposition}\label{A^*}
Let $\eta=(\eta_0,\eta_1(\cdot))\in H$. Then $\eta\in \mathcal{D}(A^*)$ if and only if $\eta_1\in W^{1,2}_{-T}$ and $\eta_1(-T)=0$. Moreover, if this is the case, then
\begin{equation}\label{astar}
A^*\eta=(r\eta_0+\eta_1(0),-\eta_1'(\cdot)).
\end{equation}
\end{proposition}
\textbf{Proof.}
Let
\[\left(\eta_0,\eta_1\right)\in\mathcal D=\left\{\eta\in H:\, \eta_1\in W_{-T}^{1,2},\,\, \eta_1(-T)=0\right\}\].
Then, for $\zeta\in \mathcal{D}(A)$,
$$\langle A\zeta,\eta\rangle=r\zeta_0\eta_0+\int_{-T}^0\zeta_1'(s)\eta_1(s)\,ds=r\zeta_0\eta_0+\zeta_0\eta_1(0)-\int_{-T}^0\zeta_1(s)\eta_1'(s)\,ds,$$
thus $\zeta\mapsto \langle A\zeta,\eta\rangle$ is continuous on $\mathcal{D}(A)$ with respect to the norm $\|\cdot\|$, i.e. $\eta\in D(A^*)$ and
$$A^*\eta=(r\eta_0+\eta_1(0),-\eta_1'(\cdot)).$$
Therefore, $\eta\in\mathcal D\left(A^*\right)$ and \eqref{astar} holds.
To show that $\mathcal D\left(A^*\right)=\mathcal D$
note first that for $t\le T$
\begin{equation}\label{sstar}
S^*(t)\left(\eta_0,\eta_1(\cdot)\right)=\left(e^{rt}\left(\eta_0+\int_{-t}^0\eta_1(\xi)e^{r\xi}d\xi\right),\eta_1(\cdot -t)I_{[-T,0]}(\cdot -t)\right).
\end{equation}
Clearly, $\mathcal D$ is dense in $H$ and it is easy to check that $S^*(t)\mathcal D\subset\mathcal D$ for any $t\ge 0$. Hence by Theorem 1.9 on p. 8 of
\cite{davies} $\mathcal D$ is dense in $\mathcal D\left(A^*\right)$ endowed with the graph norm. Finally, using \eqref{astar} it is easy to show that $\mathcal D$ is closed in the graph norm of $A^*$ and we find that $\mathcal D\left(A^*\right)=\mathcal D$.
\hfill$\square$

 \begin{lemma}\label{LipA}
 The map $F$ is Lipschitz continuous with respect to $\|\cdot\|_{-1}$.
\end{lemma}
\textbf{Proof.} Due to the Lipschitz continuity of  $f_0$, it suffices to prove that
\begin{equation}\label{iuiu}
|\eta_0|+\left|\int_{-T}^0 a(\xi) \eta_1(\xi)d\xi\right|\leq C_{a(\cdot)}\|\eta\|_{-1}, \ \ \forall\eta\in H.
\end{equation}
Indeed, since $|\eta_0|\leq r\|\eta\|_{-1}$
 $(0,a(\cdot))\in \mathcal{D}(A^*)$, we find that
\begin{eqnarray*}
\left|\int_{-T}^0 a(\xi)\eta_1(\xi)d\xi\right|&=& | \langle (0,a(\cdot)), \eta\rangle|
=| \langle (0,a(\cdot)), AA^{-1}\eta\rangle|\\
&=&| \langle A^*(0,a(\cdot)), A^{-1}\eta\rangle|
\leq \|A^*(0,a(\cdot))\|\cdot\|\eta\|_{-1}.
\end{eqnarray*}
 So, since  $|\eta_0|\leq r\|\eta\|_{-1}$,
 we get \eqref{iuiu} with $C_{a(\cdot)}=r+\|A^*(0,a(\cdot))\|$.
\hfill$\square$\\
\begin{remark}
The condition $a(-T)=0$ is in general necessary for the previous result. Indeed, consider for example the case $a(\cdot)\equiv 1$. Then the sequence
$$\eta^n=(\eta_0^n,\eta_1^n(\cdot)), \ \ \  \eta_0^n:=0, \ \eta_1^n(\cdot):=I_{[-T,-T+1/n]}(\cdot), \ \ \ n\geq 1,$$ is such that
 $$\left|\int_{-T}^0 a(\xi) \eta_1^n(\xi)d\xi\right|=1 \ \ \forall  n\geq 1, \ \ \ \ \|\eta^n\|_{-1}\rightarrow 0 \ \mbox{when} \ n\rightarrow\infty,$$
so that (\ref{iuiu}) cannot be satisfied. If for example, $f_0(r,u)=u$, the previous result does not hold.\hfill$\blacksquare$
\end{remark}
\begin{lemma}\label{estimateA}
Let $X(\cdot), \bar{X}(\cdot)$ be the mild solutions to (\ref{infinitestate}) starting respectively from $\eta, \bar{\eta}\in H$ and both under the null control. Then there exists a constant $C>0$ such that
$$\|X(t)-\bar{X}(t)\|_{-1}\leq C\|\eta-\bar{\eta}\|_{-1}, \ \ \forall t\in[0,T].$$
In particular
$$|X_0(t)-\bar{X}_0(t)|\leq rC\|\eta-\bar{\eta}\|_{-1}, \ \ \forall t\in[0,T].$$
\end{lemma}
\textbf{Proof.} From (\ref{mildcap3}) we can write, for all $t\in[0,T]$,
\begin{equation*}
X(t)-\bar{X}(t)=S(t)(\eta-\bar{\eta})+\int_0^t S(t-\tau)\left[F(X(\tau))-F(\bar{X}(\tau))\right]d\tau,
\end{equation*}
so that
\begin{equation*}
A^{-1}(X(t)-\bar{X}(t))=S(t)A^{-1}(\eta-\bar{\eta})+\int_0^t S(t-\tau)\, A^{-1}\left[F(X(\tau))-F(\bar{X}(\tau))\right]d\tau,
\end{equation*}
i.e., taking into account Lemma \ref{LipA}, there exists some $K>0$ such that
\begin{equation*}
\|X(t)-\bar{X}(t)\|_{-1}\leq K\left(\|\eta-\bar{\eta}\|_{-1}+\int_0^t \|X(\tau)-\bar{X}(\tau)\|_{-1}d\tau\right)
\end{equation*}
and the claim follows by Gronwall's Lemma.
\hfill$\square$\\

\begin{proposition}\label{usefcap3}
The set $\mathcal{D}(V)$ is open in the space $\left(H,\|\cdot\|_{-1}\right)$.
\end{proposition}
\textbf{Proof.} Let $\bar{\eta}\in \mathcal{D}(V)$, $\eta\in H_+$ and set $\bar{X}(\cdot):=X(\cdot; \bar{\eta}, 0)$, ${X}(\cdot):=X(\cdot; {\eta}, 0)$. By Proposition \ref{piu} we have $\bar{X}(t)\geq \xi>0$ for $t\in [0,T]$.
For any $\varepsilon\in\left(0,\frac{\xi}{2rC}\right)$ and any $\eta$ such that $\|\eta-\bar{\eta}\|_{-1}<\varepsilon$, Lemma \ref{estimateA} yields   $X_0(t)\geq \xi/2$ for $t\in[0,T]$. Arguing as in Proposition \ref{piu} we get $X_0(t)\geq \frac{\xi}{2}e^{-K(t-T)}$ for $t\geq T$. Thus we have the claim.\hfill$\square$
\begin{remark}
Note that $\mathcal{D}(V)$ is open also with respect to $\|\cdot\|$.
\end{remark}
\begin{proposition} \label{A^-1}
The value function is continuous with respect to $\|\cdot\|_{-1}$ on $\mathcal{D}(V)$. Moreover
\begin{equation}\label{debcont}
(\eta_n)\subset\mathcal{D}(V), \ \  \eta_n\rightharpoonup \eta\in\mathcal{D}(V) \ \ \Longrightarrow \ \ V(\eta_n)\rightarrow V(\eta).
\end{equation}
\end{proposition}
\textbf{Proof.}  The function $V$ is concave and, thanks to the proof of Lemma \ref{usefcap3}, it is $\|\cdot\|_{-1}$-locally bounded from below at the points of $\mathcal{D}(V)$. Therefore the first claim follows by a classic result of convex analysis (see e.g. \cite{Ekeland}, Chapter 1, Corollary 2.4).

 The claim \eqref{debcont} follows by  the first claim and since $A^{-1}$ is compact.\hfill$\square$\\


\subsection{Properties of superdifferential}
\label{sec:super}

In this subsection we focus on the properties of the  superdifferential of concave and $\|\cdot\|_{-1}$-continuous functions. This will be very useful in proving a regularity result for the value function. 
Recall that, if $v$ is a function defined on some open set $\mathcal{O}$ of $H$, the subdifferential and the superdifferential of $v$ at a point $\bar{\eta}\in \mathcal{O}$ are the convex and closed sets defined respectively by
$$D^-v(\bar{\eta}):=\left\{ \zeta\in H \ \Big| \ \displaystyle{\liminf_{\eta\rightarrow \bar{\eta}} \frac{v(\eta)-v(\bar{\eta})-\langle \eta-\bar{\eta}, \zeta\rangle}{\|\eta-\bar{\eta}\|}\geq 0}\right\},$$
$$D^+v(\bar{\eta}):=\left\{ \zeta\in H \ \Big| \ \displaystyle{\limsup_{\eta\rightarrow \bar{\eta}} \frac{v(\eta)-v(\bar{\eta})-\langle \eta-\bar{\eta}, \zeta\rangle}{\|\eta-\bar{\eta}\|}\leq 0}\right\}.$$
It is well-known that, if $D^+v(\eta)\cap D^-v(\eta)\neq \emptyset$, then $D^+v(\eta)\cap D^-v(\eta)=\{\zeta\}$, $v$ is differentiable at $\eta$ and $\nabla v(\eta)=\zeta$.
Moreover the set of the "reachable gradients" is defined as
$$D^*v(\bar{\eta}):=\left\{ \zeta\in H \ \Big| \ \exists \eta_n\rightarrow \bar{\eta} \ \mbox{such that} \ \exists \nabla v(\eta_n), \ \nabla v(\eta_n)\rightarrow\zeta\right\}.$$
If $\mathcal{O}$ is convex and open and $v:\mathcal{O}\rightarrow \mathbb{R}$ is concave, then the set $D^+v$ is  not empty at any point of $\mathcal{O}$ and
\begin{equation}\label{D^*}
D^+v(\bar{\eta})=\left\{ \zeta \in H \ \Big| \ v(\eta)-v(\bar{\eta})\leq \langle \eta-\bar{\eta}, \zeta\rangle, \ \ \forall\eta\in \mathcal{O}\right\}=\overline{co}\left(D^* v(\bar{\eta})\right).
\end{equation}
Moreover in this case, if $D^+v(\bar{\eta})=\{\zeta\}$, then $v$ is differentiable at $\eta$ and $\nabla v(\eta)=\zeta$.
\begin{lemma}\label{topological}
 The following statements hold:
\begin{enumerate}
\item
$A^{-1}(\mathcal{D}(V))$ is a convex open set of $(\mathcal{D}(A),\|\cdot\|)$.
\item $\mathcal{O}:= Int_{(H,\|\cdot\|)}\left(Clos_{(H,\|\cdot\|)}\left({A}^{-1}(\mathcal{D}(V))\right)\right)$ is a convex open set of $(H,\|\cdot\|)$.
\item $\mathcal{O}\supset A^{-1}(\mathcal{D}(V))$
and $\mathcal{D}(V)=\mathcal{O}\cap \mathcal{D}(A)$.
\end{enumerate}
\end{lemma}
\textbf{Proof.} The first and the second statement are obvious. We prove the third one.
Of course, since $A^{-1}(\mathcal{D}(V))$ is open in $(\mathcal{D}(A),\|\cdot\|)$, we can find $(\varepsilon_x)_{x\in A^{-1}(\mathcal{D}(V))}$, $\varepsilon_x>0$, such that
$$
A^{-1}(\mathcal{D}(V))=\bigcup_{x\in A^{-1}(\mathcal{D}(V))}B_{(\mathcal{D}(A),\|\cdot\|)}(x,\varepsilon_x).
$$
By this representation of $A^{-1}(\mathcal{D}(V))$ we can see that
$$
\mathcal{O}=\bigcup_{x\in A^{-1}(\mathcal{D}(V))}B_{(H,\|\cdot\|)}(x,\varepsilon_x).
$$
Therefore we get both the claims of the third statement.\hfill$\square$

\begin{proposition}\label{importante}
 Let $v:\mathcal{D}(V)\rightarrow \mathbb{R}$ be a concave function continuous with respect to $\|\cdot\|_{-1}$. Then
\begin{enumerate}
\item $v=u\circ A^{-1}$, where
$u:\mathcal{O}\subset H \rightarrow \mathbb{R}$
is a concave $\|\cdot\|$-continuous function.
\item $D^+v(\eta)\subset \mathcal{D}(A^*)$, for any $\eta\in \mathcal{D}(V)$.
\item $D^+ u(A^{-1}\eta)= A^*D^+v(\eta)$, for any $\eta\in\mathcal{D}(V)$. In particular, since $A^*$ is injective, $v$ is differentiable at $\eta$ if and only if $u$ is differentiable at $A^{-1}\eta$.
\item If $\zeta\in D^*v(\eta)$, then there exists a sequence $\eta_n\rightarrow \eta$ such that there exist $\nabla v(\eta_n)$, $\nabla v(\eta_n)\rightarrow \zeta$ and $A^*\nabla v(\eta_n)\rightharpoonup A^*\zeta$.
\end{enumerate}
\end{proposition}

\textbf{Proof.} Within this proof, for $\eta\in\mathcal{D}(V)$, we set $\eta':=A^{-1}\eta$. Since $A^{-1}$ is one-to-one, there is a one-to-one correspondence between the elements $\eta\in\mathcal{D}(V)$ and $\eta'\in A^{-1}(\mathcal{D}(V))$.

\textbf{1.} Let us define the function $u_0:A^{-1}(\mathcal{D}(V))\rightarrow \mathbb{R}$ by
$$u_0(\eta'):=v(\eta).$$
Thanks to the assumptions on $v$, $u_0$ is a concave continuous function on $(A^{-1}(\mathcal{D}(V)),\|\cdot\|)$. By the third statement of Lemma \ref{topological} we see that $A^{-1}(\mathcal{D}(V))$ is $\|\cdot\|$-dense in $\mathcal{O}$. Since $v$ is concave it is locally Lipschitz continuous, so that can be extended to a concave $\|\cdot\|$-continuous function $u$ defined on $\mathcal{O}$. This function $u$ satisfies the claim by construction.

\textbf{2.} Let $\bar{\eta}\in \mathcal{D}(V)$, $\zeta\in D^+v(\bar{\eta})$. Then
$$
v(\eta)-v(\bar{\eta})\leq \langle  \eta-\bar{\eta}, \zeta \rangle, \ \ \ \ \ \forall \eta\in \mathcal{D}(V),
$$
i.e.
$$
u(\eta')-u(\bar{\eta}')\leq \langle  A(\eta'-\bar{\eta}'), \zeta \rangle, \ \ \ \ \ \forall \eta'\in A^{-1}(\mathcal{D}(V)).
$$
Thus the function
$$
\begin{array}{cccc}
T_\zeta:&(\mathcal{D}(A),\|\cdot\|)&\longrightarrow & \mathbb{R},\\
&\eta'&\longmapsto &\langle  A\eta',\zeta\rangle,
\end{array}
$$
is lower semicontinuous at $\bar{\eta}'$. It is also linear and therefore it is continuous on $(\mathcal{D}(A), \|\cdot\|)$, so that we can conclude that $\zeta\in \mathcal{D}(A^*)$.

\textbf{3.} Let $\bar{\eta}\in \mathcal{D}(V)$, $\zeta\in D^+v(\bar{\eta})$. Then
$$
v(\eta)-v(\bar{\eta})\leq \langle  \eta-\bar{\eta}, \zeta \rangle, \ \ \ \ \ \forall \eta\in \mathcal{D}(V),
$$
i.e.
$$
u(\eta')-u(\bar{\eta}')\leq \langle  A(\eta'-\bar{\eta}'), \zeta \rangle= \langle  (\eta'-\bar{\eta}'), A^*\zeta \rangle, \ \ \ \ \ \forall \eta'\in A^{-1}(\mathcal{D}(V)),
$$
so that $A^*\zeta\in D^+u(\bar{\eta}')$, which gives $D^+ u(A^{-1}\eta)\supset A^*D^+v(\eta)$.

Conversely let $\bar{\eta}'\in \mathcal{A}^{-1}(\mathcal{D}(V))$ and $\zeta'\in D^+u(\bar{\eta}')$. Then
$$
u(\eta')-u(\bar{\eta}')\leq \langle  A(\eta'-\bar{\eta}'), \zeta' \rangle, \ \ \  \forall \eta'\in A^{-1}(\mathcal{D}(V)),
$$
i.e.
$$
v(\eta)-v(\bar{\eta})\leq \langle  A^{-1}(\eta-\bar{\eta}), \zeta' \rangle=  \langle(\eta-\bar{\eta}), (A^{-1})^*\zeta' \rangle, \ \ \ \ \ \forall \eta\in \mathcal{D}(V).
$$
Since $(A^{-1})^*=(A^*)^{-1}$, we get $(A^*)^{-1}\zeta'\in D^+v(\bar{\eta})$, which gives $D^+ u(A^{-1}\eta)\subset A^*D^+v(\eta)$.

\textbf{4.} Let $\bar{\eta}\in \mathcal{D}(V)$ and $\zeta\in D^*v(\bar{\eta})$. Due to (\ref{D^*}), we can find a sequence $\eta_n\rightarrow\bar{\eta}$ such that $\nabla v(\eta_n)$ exists for any $n\in\mathbb{N}$ and $\nabla v(\eta_n)\rightarrow \zeta$. Thanks to the third claim we can say that also $\nabla u(\bar{\eta}'_n)$ exists and $\nabla u(\bar{\eta}_n')=A^*\nabla v(\eta_n)$. The sequence $\nabla u(\bar{\eta}_n')$ is bounded, due to the fact that the set-valued map $\eta'\mapsto D^+u(\bar{\eta}')$ is locally bounded. Therefore from any subsequence we can extract a subsubsequence weakly converging to some element $\zeta'\in H$. $A^*$ is a closed operator, so that it is also a weakly closed operator. Therefore we can conclude that $\zeta\in \mathcal{D}(A^*)$ and $\zeta'=A^*\zeta$. Since this holds for any subsequence, we can conclude that $A^*\nabla v(\eta_n)\rightharpoonup A^*\zeta$.
\hfill$\square$

\section{Dynamic Programming}

The dynamic programming principle states that, for any $\eta\in
\mathcal{D}(V)$ and for any $s\geq 0$,
\begin{equation*}
V(\eta)=\sup_{c(\cdot)\in \mathcal{C}(\eta)}\left[\int_{0}^se^{-\rho t}\left(U_1(c(t)+U_2(X_0(t))\right)dt+e^{-\rho s}V(X(s;\eta,c(\cdot)))\right].
\end{equation*}
Its differential version is the Hamilton-Jacobi-Bellman (from now on HJB) equation on $\mathcal{D}(V)$, which in our case reads as
\begin{equation}\label{HJB}
\rho v(\eta)= \langle \eta, A^* \nabla v(\eta)\rangle +f(\eta)v_{\eta_0}(\eta)+ U_2(\eta_0)+\mathcal{H}(v_{\eta_0}(\eta)),
\end{equation}
where $\mathcal{H}$ is the Legendre transform of $U_1$, i.e.
$$\mathcal{H}(\zeta_0):=\sup_{c\geq 0}\left(U_1(c)-\zeta_0 c\right), \  \ \ \zeta_0> 0.$$
Due to  Hyphothesis \ref{hyp:utility}-(i) and to Corollary 26.4.1 of \cite{Rock}, we have that $\mathcal{H}$ is strictly convex on $(0,+\infty)$. Notice that, thanks to Proposition \ref{strict}-(3),
$$D^+_{\eta_0}V(\eta):=\{\zeta_0 \in \mathbb{R} \ | \ (\zeta_0,\zeta_1(\cdot))\in D^+V(\eta)\}\subset (0,\infty)$$
for any $\eta\in\mathcal{D}(V)$, i.e. where $\mathcal{H}$ is defined.

\subsection{Viscosity solutions}\label{sec:visc}
First we study the HJB equation using the viscosity solutions
approach. In order to follow this approach, we have to define a
suitable set of regular test functions. This is the set
\begin{equation}\label{testcap3}
\tau:=\Big\{\varphi\in C^1(H) \ | \ \nabla \varphi (\cdot)\in \mathcal{D}(A^*), \ \eta_n\rightarrow \eta\Rightarrow A^*\nabla\varphi(\eta_n)\rightharpoonup A^*\nabla\varphi(\eta)\Big\}.
\end{equation}
Let us define, for $c\geq 0$,  the operator $\mathcal{L}^{c}$ on $\tau$ by
$$[\mathcal{L}^{c}\varphi](\eta):=-\rho \varphi (\eta)+\langle \eta ,A^*\nabla\varphi(\eta)\rangle+  f(\eta)\varphi_{\eta_0}(\eta)-c\varphi_{\eta_0}(\eta).$$

\begin{lemma}\label{identity}
Let $\varphi\in\tau$, $c(\cdot)\in L^1([0,+\infty);\mathbb{R}^+)$ and set $X(t):=X(t;\eta,c(\cdot))$. Then the following identity holds for any $t\geq 0$:
$$e^{-\rho t}\varphi(X(t))-\varphi(\eta)=\int_0^te^{-\rho s}[\mathcal{L}^{c(s)}\varphi](X(s))ds.$$
\end{lemma}

\textbf{Proof.} The statement holds if we replace $A$ with the Yosida approximations. Then we can pass to the limit and get the claim thanks to the regularity properties of the functions belonging to  $\tau$.\hfill$\square$

\begin{definition}\label{def:viscosity}

 (i) A continuous function $v:\mathcal{D}(V)\rightarrow \mathbb{R}$ is called a viscosity subsolution of (\ref{HJB}) on $\mathcal{D}(V)$ if for any $\varphi\in\tau$ and any $\eta_M\in\mathcal{D}(V)$ such that $v-\varphi$ has a $\|\cdot\|$-local maximum at $\eta_M$ we have
$$\rho v(\eta_M)\leq   \langle \eta_M, A^* \nabla \varphi(\eta_M)\rangle +f(\eta_M)\varphi_{\eta_0}(\eta_M)+U_2(\eta_0)+ \mathcal{H}(\varphi_{\eta_0}(\eta_M)).$$

(ii) A continuous function $v:\mathcal{D}(V)\rightarrow \mathbb{R}$ is called a viscosity supersolution of (\ref{HJB}) on $\mathcal{D}(V)$ if for any $\varphi\in\tau$ and any $\eta_m\in\mathcal{D}(V)$ such that $v-\varphi$ has a $\|\cdot\|$-local minimum at $\eta_m$  we have
$$\rho v(\eta_m)\geq   \langle \eta_m, A^* \nabla \varphi(\eta_m)\rangle +f(\eta_m)\varphi_{\eta_0}(\eta_m)+ U_2(\eta_0)+ \mathcal{H}(\varphi_{\eta_0}(\eta_m)).$$

(iii) A continuous function $v:\mathcal{D}(V)\rightarrow \mathbb{R}$ is called a viscosity supersolution of (\ref{HJB}) on $\mathcal{D}(V)$ if it is both a viscosity sub and supersolution.

\end{definition}

We can prove the following:
\begin{theorem}\label{TH:visc}
The value function $V$ is a viscosity solution of (\ref{HJB}) on $\mathcal{D}(V)$.
\end{theorem}
\textbf{Proof.} (i) We prove that $V$ is a viscosity subsolution. Let $(\eta_M,\varphi)\in \mathcal{D}(V)\times \tau$ be such that $V-\varphi$ has a local maximum at $\eta_M$. Without loss of generality we can suppose $V(\eta_M)=\varphi(\eta_M)$. Let us suppose, by contradiction that there exists $\nu>0$ such that
$$2\nu\leq \rho V(\eta_M) -\big( \langle \eta_M, A^* \nabla \varphi(\eta_M)\rangle +f(\eta_M)\varphi_{\eta_0}(\eta_M)+U_2(\eta_0)+ \mathcal{H}(\varphi_{\eta_0}(\eta_M))\big).$$
Let us define the function
$$\tilde{\varphi}(\eta):= V(\eta_M)+\langle \nabla \varphi(\eta_M),\eta-\eta_M\rangle+\|\eta-\eta_M\|_{-1}^2.$$
We have
$$\nabla\tilde{\varphi}(\eta)=\nabla\varphi(\eta_M)+(A^*)^{-1}A^{-1}(\eta-{\eta_M}), $$
Thus $\tilde{\varphi}$ is a test function and we must have also
$$2\nu\leq \rho V(\eta_M) -\big( \langle \eta_M, A^* \nabla \tilde{\varphi}(\eta_M)\rangle +f(\eta_M)\tilde{\varphi}_{\eta_0}(\eta_M)+U_2(\eta_0)+ \mathcal{H}(\tilde{\varphi}_{\eta_0}(\eta_M))\big).$$
By concavity of $V$ we have
$$V(\eta_M)=\tilde{\varphi}(\eta_M), \ \ \ \ \ \ \tilde{\varphi}(\eta)\geq  V(\eta)+\|\eta-\eta_M\|_{-1}^2,\ \eta\in\mathcal{D}(V).$$
By the continuity property of $\tilde{\varphi}$ we can find $\varepsilon>0$ such that
$$\nu\leq \rho V(\eta) -\big( \langle \eta, A^* \nabla \tilde{\varphi}(\eta)\rangle +f(\eta)\tilde{\varphi}_{\eta_0}(\eta)+U_2(\eta_0)+ \mathcal{H}(\tilde{\varphi}_{\eta_0}(\eta))\big), \ \ \ \ \eta\in B(\eta_M,\varepsilon).$$
Take a sequence $\delta_n>0$, $\delta_n\rightarrow 0$ and, for any $n$, take a $\delta_n$-optimal control $c_n
(\cdot)\in \mathcal{C}_{ad}(\eta_M)$. Set $X^n(\cdot):=X(\cdot;\eta_M, c_n
(\cdot))$ and define
$$t_n
:=\inf\{t\geq 0 \ | \ \|X^n(t)-\eta_M\|=\varepsilon\}\wedge 1.$$
Of course $t_n
$ is well-defined and belongs to $(0,1]$. Moreover, by continuity of trajectories,  $X^n(t)\in B(\eta_M,\varepsilon)$, for $t\in [0,t_n
)$.
We distinguish two cases:
$$\limsup_n t_n
=0, \ \ \ \ \ \ \mbox{or} \ \ \ \ \ \ \limsup_n t_n
>0.$$
In the first case we can write
\begin{eqnarray*}
\delta_n &\geq &-\int_0^{t_n
} e^{-\rho t} \left[ U_1(c_n
(t))+ U_2(X^n_0(t))\right]dt-\left(e^{-\rho t_n
}V(X(t_n))-V(\eta_M)\right)\\
&\geq & -\int_0^{t_n
} e^{-\rho t} \left[ U_1(c_n
(t))+ U_2(X^n_0(t))\right]dt-\left(e^{-\rho t_n
}(\tilde{\varphi}(X^n(t_n
)))-\tilde{\varphi}(\eta_M)\right)+ e^{-\rho t_n
}\|X^n(t_n
)-\eta_M\|_{-1}^2\\
&= & -\int_0^{t_n
} e^{-\rho t} \left[ U_1(c_n
(t))+ U_2(X^n_0(t))+[\mathcal{L}^{c_n
(t)}\tilde{\varphi}] (X^n(t)) \right]dt + e^{-\rho t_n
}\|X^n(t_n
)-\eta_M\|_{-1}^2\\
&\geq & -\int_0^{t_n
} e^{-\rho t} \Big[U_2(X^n_0(t))-\rho \tilde{\varphi} (X^n(t))+\langle X^n(t), A^*\nabla\tilde{\varphi}(X^n(t))\rangle\\
&& \ \ \ \ \ \ \ \ \ +  f(X^n(t))\tilde{\varphi}_{\eta_0}(X^n(t))+\mathcal{H}(\tilde{\varphi}_{\eta_0}(X^n(t))) \Big]dt+ e^{-\rho t_n
}\|X^n(t_n
)-\eta_M\|_{-1}^2\\
&\geq& t_n
\nu + e^{-\rho t_n
}\|X^n(t_n
)-\eta_M\|_{-1}^2,
\end{eqnarray*}
thus it has to be
$$\|X^n(t_n
)-\eta_M\|_{-1}^2\rightarrow 0.$$
Let us show that this is impossible. The above convergence implies in particular that
\begin{equation}\label{conv1}
|X^n_0(t_n
)-(\eta_M)_0|\rightarrow 0.
\end{equation}
Moreover, by definition of $t_n$, it has to be
\begin{equation}\label{dom}
|X^n_0(t
)-(\eta_M)_0|\leq \varepsilon, \  \ \ t\in[0,t_n].
\end{equation}
Since $t_n
\rightarrow 0$, taking into account (\ref{dom}), we have also
\begin{equation}\label{conv2}
\|X^n_1(t_n
)-(\eta_M)_1\|_{L^2_{-T}}\rightarrow 0.
\end{equation}
The convergences (\ref{conv1}) and (\ref{conv2}) are not compatible with the definition of $t_n$ and the contradiction arises.
In the second case we can suppose, eventually passing to a subsequence, that $t_n
\rightarrow \bar{t}\in(0,1]$. So we get as before

\begin{equation*}
\delta_n  \geq t_n
\nu +e^{-\rho t_n}\|X^n(t_n)-\eta_M\|_{-1}^2\geq t_n\nu;
\end{equation*}
since $\delta_n\rightarrow 0$ and $t_n\nu \rightarrow \bar{t}\nu$, again a contradiction arises.
\smallskip

(ii) The proof that $V$ is a viscosity supersolution is standard, see e.g. \cite{LiYong}. \hfill$\square$

\subsection{Smoothness of viscosity solutions}\label{sec:reg}
In this subsection we show that the concave $\|\cdot\|_{-1}$-continuous viscosity solutions of (\ref{HJB}) (so that in particular the value function $V$) are differentiable along the direction $\hat{n}=(1,0)$. For this purpose  we need the following lemma.

\begin{lemma}\label{visc=}
Let $v:\mathcal{D}(V)\rightarrow\mathbb{R}$ a concave $\|\cdot\|_{-1}$-continuous function and suppose that $\bar{\eta}\in \mathcal{D}(V)$ is a differentiability point for $v$ and that $\nabla v(\bar{\eta})=\zeta$. Then
\begin{enumerate}
 \item There exists a test function $\varphi$ such that $v-\varphi$ has a local maximum at $\bar{\eta}$ and $\nabla \varphi(\bar{\eta})=\zeta$.
 \item There exists a test function $\varphi$ such that $v-\varphi$ has a local minimum at $\bar{\eta}$ and $\nabla \varphi(\bar{\eta})=\zeta$.
\end{enumerate}
\end{lemma}
\textbf{Proof.} Thanks to Proposition \ref{importante} and due to the concavity of $v$, the first statement is clearly satisfied by the function $\langle \cdot,\zeta\rangle$. We prove now the second statement, which is more delicate. We use the notation of Proposition \ref{importante}. Thanks to the third claim of  Proposition \ref{importante}, we have $A^*\zeta\in D^+ u(\bar{\eta}')$. This means that
$$u(\eta')-u(\bar{\eta}')-\langle \eta'-\bar{\eta}', A^*\zeta\rangle\geq -\|\eta'-\bar{\eta}'\| \cdot\varepsilon (\|\eta'-\bar{\eta}'\|),$$
where $\varepsilon: [0,+\infty)\rightarrow [0,+\infty)$ is an increasing function such that $\varepsilon (\|\eta'-\bar{\eta}'\|)\rightarrow 0$, when $\|\eta'-\bar{\eta}'\|\rightarrow 0$.
The previous inequality can be rewritten also as
$$u(\eta')-u(\bar{\eta}')-\langle A(\eta'-\bar{\eta}'), \zeta\rangle\geq -\|\eta'-\bar{\eta}'\|\cdot \varepsilon (\|\eta'-\bar{\eta}'\|).$$
Passing to $v$ this reads as
$$v(\eta)-v(\bar{\eta})-\langle \eta-\bar{\eta}, \zeta\rangle\geq -\|\eta-\bar{\eta}\|_{-1} \cdot\varepsilon \left(\|\eta-\bar{\eta}\|_{-1}\right),$$
where $\varepsilon \left(\|\eta-\bar{\eta}\|_{-1}\right)\rightarrow 0$, when $\|\eta-\bar{\eta}\|_{-1}\rightarrow 0$. We look for a test function of this form:
$$\varphi(\eta)=v(\bar{\eta})+\langle \eta-\bar{\eta}, \zeta\rangle + g\left(\|\eta-\bar{\eta}\|_{-1}\right),$$
where $g:[0,+\infty)\rightarrow[0,+\infty)$ is a suitable increasing $C^1$ function such that $g(0)=0$. Notice that $\varphi (\bar{\eta})=v(\bar{\eta})$, so that, in order to prove that $v-\varphi$ has a local minimum at $\bar{\eta}$, we have to prove that $\varphi\leq v$ in a neighborhood of $\bar{\eta}$.

Let us define the function
$$g(r):=\int_0^{2r}\varepsilon (s)ds\geq\int_r^{2r}\varepsilon(s)ds\geq r\varepsilon(r).$$
Then
\begin{eqnarray*}
\varphi(\eta)&=&v(\bar{\eta}+\langle\eta-\bar{\eta},\zeta\rangle- g(\|\eta-\bar{\eta}\|_{-1})\\
&\leq &v(\bar{\eta})+\langle\eta-\bar{\eta},\zeta\rangle-\|\eta-\bar{\eta}\|_{A{-1}} \cdot\varepsilon(\|\eta-\bar{\eta}\|_{-1})\\
&\leq & v(\eta).
\end{eqnarray*}
Moreover $$\nabla\varphi(\eta)=
\begin{cases}
\displaystyle{\zeta- (A^*)^{-1}\frac{\varepsilon (2\|\eta-\bar{\eta}\|_{-1})}{\|\eta-\bar{\eta}\|_{-1}}A^{-1}(\eta-\bar{\eta}),\ \ \ \ \mbox{if} \ \eta\neq \bar{\eta},}\\
\zeta, \qquad\qquad\qquad\qquad\qquad\qquad\qquad\qquad \quad\mbox{if} \ \eta=\bar{\eta},
\end{cases}$$
i.e. $\varphi$ is a test function and $\nabla\varphi(\bar{\eta})=\zeta$.\hfill$\square$\\

Now we can state and prove the main result.
\begin{theorem}\label{TH:reg}
 Let $v$ be a concave $\|\cdot\|_{-1}$-continuous viscosity solution of (\ref{HJB}) on $\mathcal{D}(V)$. Then $v$ is differentiable along the direction $\hat{n}=(1,0)$ at any point $\eta\in\mathcal{D}(V)$ and the function $\eta\mapsto v_{\eta_0}(\eta)$ is continuous on $\mathcal{D}(V)$.
\end{theorem}
\textbf{Proof.} Let ${\eta}\in\mathcal{D}(V)$ and $\zeta, \xi\in D^*v(\bar{\eta})$. Thanks to Proposition \ref{importante}, there exist sequences $(\eta_n)$, $(\tilde{\eta}_n)$ such that:
\begin{itemize}
\item $\eta_n\rightarrow {\eta}$, $\tilde{\eta}_n\rightarrow {\eta}$;
\item $\nabla v(\eta_n)$ and $\nabla v(\tilde{\eta}_n)$ exist for all $n\in \mathbb{N}$;
\item $A^* \nabla v(\eta_n)\rightharpoonup A^*\zeta$ and $A^* \nabla v(\tilde{\eta}_n)\rightharpoonup A^*\xi$.
\end{itemize}
Thanks to Lemma \ref{visc=} we can write, for any $n\in\mathbb{N}$,
$$\rho v(\eta_n)= \langle \eta_n, A^* \nabla v(\eta_n)\rangle +f(\eta_n)v_{\eta_0}(\eta_n)+U_2(\eta^n_0) +\mathcal{H}(v_{\eta_0}(\eta_n)),$$
$$\rho v(\tilde{\eta}_n)= \langle \tilde{\eta}_n, A^* \nabla v(\tilde{\eta}_n)\rangle +f(\tilde{\eta}_n)v_{\eta_0}(\tilde{\eta}_n)+ U_2(\eta_0^n)+\mathcal{H}(v_{\eta_0}(\tilde{\eta}_n)).$$
Passing to the limit we get
\begin{equation}\label{eq1}
\langle {\eta}, A^* \zeta\rangle +f({\eta})\zeta_0+U_2({\eta}_0)+\mathcal{H}(\zeta_0)=\rho v({\eta})=  \langle {\eta}, A^* \xi\rangle +f({\eta})\xi_0+U_2({\eta}_0)+ \mathcal{H}(\xi_0).
\end{equation}
On the other hand $\lambda \zeta+(1-\lambda)\xi\in D^+ v(\bar{\eta})$, for any $\lambda\in(0,1)$, so that we have the subsolution inequality
\begin{equation}\label{eq2}
\rho v({\eta})\leq \langle  {\eta}, A^* [\lambda\zeta+(1-\lambda)\xi]\rangle +f({\eta})[\lambda\zeta_0+(1-\lambda)\xi_0]+ U_2({\eta}_0)+\mathcal{H}(\lambda\zeta_0+(1-\lambda)\xi_0), \ \ \forall\lambda\in(0,1).
\end{equation}
Combining (\ref{eq1}) and (\ref{eq2}) we get
$$\mathcal{H}(\lambda\zeta_0+(1-\lambda)\xi_0)\geq \lambda\mathcal{H}(\zeta_0)+(1-\lambda)\mathcal{H}(\xi_0);$$
since $\mathcal{H}$ is strictly convex, the previous inequality implies $\zeta_0=\xi_0$. This means that the projection of $D^* v({\eta})$ onto $\hat{n}$ is a singleton. Thanks to (\ref{D^*}) this implies also that the projection of $D^+v({\eta})$ onto $\hat{n}$ is a singleton and therefore that $v$ is differentiable in the direction $\hat{n}$ at ${\eta}$.
\smallskip

We prove now that the map $\eta\mapsto v_{\eta_0}(\eta)$ is continuous on $\mathcal{D}(V)$. Let ${\eta}\in\mathcal{D}(V)$ and let $(\eta^n)$ be a sequence such that $\eta^n\rightarrow{\eta}$. We have to show that
$v_{\eta_0}(\eta^n)\rightarrow v_{\eta_0}({\eta})$. Of course for any $n\in\mathbb{N}$ there exists $p_1^n\in L^2_{-T}$ such that $(v_{\eta_0}(\eta^n),p_1^n)\in D^+ v(\eta^n)$. Since $v$ is concave, it is also locally Lipschitz continuous so that the super-differential is locally bounded. Therefore, from any subsequence $\left(v_{\eta_0}(\eta^{n_k})\right)$, we can extract a sub-subsequence $\left(v_{\eta_0}(\eta^{n_{k_h}})\right)$ such that $\left(v_{\eta_0}(\eta^{n_{k_h}}),p_1^{n_{k_h}}\right)$ is weakly convergent towards some limit point. Due to the concavity of $v$ this limit point must live in the set $D^+ v(\eta)$. In particular the limit point of $\left(v_{\eta_0}(\eta^{n_{k_h}})\right)$ must coincide with $v_{\eta_0}({\eta})$. This holds true for any subsequence $\left(v_{\eta_0}(\eta^{n_k})\right)$, so that the claim follows by the usual argument on subsequences.  \hfill$\square$
\begin{remark} Notice that in the assumptions of Theorem \ref{TH:reg} we do not require that $v$ is the value function, but only  that it is a concave $\|\cdot\|_{-1}$-continuous viscosity solution of \eqref{HJB}.\hfill $\blacksquare$
\end{remark}

\begin{remark}
Thanks to the regularity result of the previous subsection 
 we can
define the  feedback map
\begin{equation*}\label{feedbackmap}
C(\eta):=\mbox{argmax}_{c\geq 0} \left(U_1(c)-c V_{\eta_0}(\eta)\right), \ \ \ \eta\in \mathcal{D}(V),
\end{equation*}
which, at least formally, should define an optimal strategy for the problem.

The object of the forthcoming paper \cite{FGG2} is the study of the
closed loop equation associated to this map and the proof of a
Verification Theorem showing that this map actually defines an
optimal  feedback strategy for the problem.\hfill$\blacksquare$
\end{remark}

\begin{remark}\label{rm:noBambi}
When the delay is concentrated in a point in a linear way, we could
tempted to insert the delay term in the infinitesimal generator $A$
and try to work as done in Section \ref{sec:infinite}. Unfortunately
this is not possible. Indeed consider this simple case:
\begin{equation*}
\begin{cases}
y'(t)=ry(t)+y\left(t-T\right),\\
y(0)=\eta_0, \ y(s)=\eta_1(s), \ s\in[-T,0),
\end{cases}
\end{equation*}
In this case we can define
$$
A:\mathcal{D}(A)\subset H\longrightarrow H, \qquad
(\eta_0,\eta_1(\cdot))\longmapsto  (r
\eta_0+\eta_1(-T),\eta_1'(\cdot)).
$$
where again
$$\mathcal{D}(A):=\{ \eta\in H \ | \ \eta_1(\cdot)\in W^{1,2}([-T,0]; \mathbb{R}), \ \eta_1(0)=\eta_0\}.$$
The inverse of $A$ is the operator
$$
A^{-1}:(H,\|\cdot\|)\longrightarrow (\mathcal{D}(A),\|\cdot\|)
\qquad(\eta_0,\eta_1(\cdot))\longmapsto
\displaystyle{\left(\frac{\eta_0-c}{r}, \ c+\int_{-T}^\cdot
\eta_1(\xi)d\xi\right),}
$$
where
$$c=\frac{1}{r+1}\,\eta_0-\frac{r}{r+1}\,\int_{-T}^0\eta_1(\xi)d\xi.$$
In this case we would have the first part of Lemma \ref{estimateA},
but not the second part, because it is not possible to control
$|\eta_0|$ by $\|\eta\|_{-1}$. Indeed take for example $r$ such that
$\frac{1-r}{1+r}=\frac{1}{2}$, and $(\eta^n)_{n\in\mathbb{N}}\subset
H$ such that
$$\eta_0^n=1/2, \ \ \ \int_{-T}^0\eta_1^n(\xi)d\xi=1, \ \ \ n\in\mathbb{N}.$$
We would have $c=1/2$, so that $\left|\frac{\eta_0^n-c}{r}\right|=0$. Moreover we can choose $\eta_1^n$ such that, when $n\rightarrow \infty$,
$$\int_{-T}^0\left|\,\frac{1}{2}+\int_{-T}^s\eta_1^n(\xi)d\xi\right|^2ds\longrightarrow 0.$$
Therefore we would have $|\eta_0^n|=1/2$ and $\|\eta^n\|_{-1}\rightarrow 0$. This shows that the second part of Lemma \ref{estimateA} does not hold. Once this part does not hold, then everything in the following argument breaks down.\hfill$\blacksquare$
\end{remark}


\end{document}